\documentclass[10pt, oneside]{article}
\usepackage{lineno}
\usepackage[a4paper, margin=1in]{geometry}
\usepackage{amsmath}
\usepackage{amssymb,amsfonts}
\usepackage{cool}
\usepackage{psfrag}
\usepackage{enumerate}
\usepackage{subfigure}
\usepackage{epsfig}
\usepackage{color}
\usepackage[usenames,dvipsnames]{xcolor}
\newcommand{\mat}[1]{\mbox{\boldmath{$#1$}}} 
\usepackage{lineno,hyperref}
\usepackage{float}
\modulolinenumbers[5]

\usepackage{natbib}
\hypersetup{colorlinks,urlcolor=blue,citecolor=blue,linkcolor=blue} 

\usepackage{booktabs}

\usepackage{fullpage}
\usepackage{setspace}
\usepackage{flushend}
\usepackage{multicol}
\usepackage{multirow}
\usepackage{rotating}
\usepackage[lined,linesnumberedhidden]{algorithm2e}

\usepackage{url}
\usepackage[normalem]{ulem}
\usepackage{bm}

\usepackage[labelsep=period]{caption}

\usepackage{authblk}
\usepackage{epstopdf}

\begin{document}



\title{Bartlett and Bartlett-type corrections in heteroscedastic symmetric nonlinear regression models}

\author{
Mariana C. Ara\'ujo\textsuperscript{1}
\thanks{
Correspondence to: Mariana C. Ara\'ujo \\ Adress: Departamento de Estat\'istica, UFRN, CCET, Lagoa Nova, Av. Senador Salgado Filho, 3000 , 59078-970, Natal, RN , Brazil  \\ E-mail: mariana@ccet.ufrn.br \\ Phone number: +558433422264}
\ , Audrey H.M.A. Cysneiros\textsuperscript{2}
\ , Lourdes C. Montenegro\textsuperscript{3}}
\affil{\textsuperscript{1}Departamento de Estat\'istica, Universidade Federal do Rio Grande do Norte, Lagoa Nova, Av. Senador Salgado Filho, 3000 , 59078-970, Natal, RN , Brazil \\ \textsuperscript{2}Departamento de Estat\'istica, Universidade Federal de Pernambuco, Cidade Universitária, Av. Prof. Moraes Rego, 1235, 50740-540, Recife, PE, Brazil \\ \textsuperscript{3}Departamento de Estat\'istica, Universidade Federal de Minas Gerais, Pampulha, Av. Pres. Ant\^onio Carlos, 6627, 31270-901, Belo Horizonte, MG, Brazil  }
\date{}

\maketitle

\begin{abstract}
\hspace{-.5cm}This paper provides general expression for Bartlett and Bartlett-type correction factors for the likelihood ratio and gradient statistics to test the dispersion parameter in heteroscedastic symmetric nonlinear models. This class of regression models is potentially useful for modeling  data containing outlying observations. We consider a partition on the dispersion parameter vector in order to test the parameters of interest. Furthermore, we develop Monte Carlo simulations to compare the finite sample performances of the corrected tests proposed with the usual and modified score tests, likelihood and gradient tests, the Bartlett-type corrected score test and bootstrap corrected tests. Our simulation results favor the score and gradient corrected tests as well as the bootstrap tests. An empirical application is presented for illustrative purposes.
\vspace{0.5cm}

\hspace{-.5cm}\textit{keywords:} Bartlett corrections, Bartlett-type corrections, Bootstrap, Gradient test, Large-sample test statistics.
\end{abstract}


\section{Introduction}
The symmetric class of models has received increasing attention in the literature. Including the normal distribution, the symmetric family covers both light and heavy tailed distributions including Cauchy, Student$-t$, generalized Student$-t$ and power exponential, among others. The symmetric models provide a very useful extension of the normal model, once that using a heavy tailed distribution for the error component reduces the influence of extreme observations and enables carrying out a more robust statistical analysis (\citealt{llt1989}). An extensive range of practical applications considering symmetric distributions can be found in various fields, such as engineering, biology and economics, among others. The symmetric regression models have been being subject of several studies (e.g., \citealt{lzx2009}; \citealt{ccc2010}; \citealt{l2012}; \citealt{mc2018}). 

Constant dispersion is often a standard assumption when symmetric data are fitted. However, in many practical situations this condition is not satisfied, requiring verification, since the inference strategies change when one observes variable dispersion over the observations. The likelihood ratio (LR), Wald and score are the large-sample tests commonly used for this purpose. The recently proposed gradient test (\citealt{t2002}), whose statistic shares the same first order asymptotic properties with the LR, Wald and score statistics  (\citealt{lf2012a}), has been the subject of many studies in the past few years (e.g., \citealt{l2011}; \citealt{l2013}; \citealt{lf2012b}; \citealt{mf2017}), given that when compared to the Wald and score statistics, the gradient statistic does not depend on the information matrix, either expected or observed, and is also simpler to compute.

The four statistics for testing hypothesis in regression models have the null asymptotic $\chi^2_q$ distribution, where $q$ is the difference between the dimensions of the parameter space under the two hypotheses being tested, up to an error order $n^{-1}$. Relying on inference in tests based on such statistics has less justification when dealing with small and moderate sized samples. A strategy to improve the $\chi^2$ approximation for the exact distributions of the LR, score and gradient statistics is to multiply them by a correction factor. For the LR statistic, \citealt{b1937} proposed a correction factor known as the Bartlett correction, which was put into a general framework later by \citealt{l1956}, while for the score statistic, \citealt{cf1991} proposed a Bartlett-type correction. Based on the results from \citealt{cf1991}, a Bartlett-type correction of the gradient statistic was recently proposed in a general framework by \citealt{vfl2013}. For the Wald statistic, there is no Bartlett or Bartlett-type correction to improve the approximation from its exact distribution to the $\chi^2$ distribution in a general setting. The corrected versions of the test statistics have the same $\chi^2_q$ null distribution with approximation error of order $n^{-2}.$ \citealt{gc2014} shows additional details on Bartlett corrections. Improved tests have been discussed in some recent articles, in particular \citealt{lcm2012}, \citealt{bc2013}, \citealt{vfl2014} and \citealt{mfl2017}.


Considering the class of heteroscedastic symmetric nonlinear models {\footnotesize (HSNLM)} proposed by \citealt{ccc2010}, \citealt{n2010} derived a Bartlett-type correction for the score statistic, proceeding with a numerical study to test the regression coefficients in the dispersion parameter. In this paper, our main goal is to derive Bartlett and Bartlett-type corrections to improve inference on the dispersion parameter based on the LR and gradient statistics, respectively, for the class of HSNLM considering the parametrization presented in \citealt{ccc2010}. Furthermore, we consider a partition of the dispersion parameter which is an advantage, since that in some cases we are not interest in making inference on all parameters of the model. It is important to mention that one of the main results presented in this paper, which is the Bartlett correction factor for the $LR$ statistic, is not the same presented in \citealt{acm2020}.

In order to achieve our aim, we adopt a regression structure to model the dispersion parameter vector so that under the null hypothesis the dispersion is constant. In other words, the null hypothesis delivers the symmetric nonlinear regression model. Our results provide a new class of tests which can be used in practical applications, mainly those involving small datasets. 

We perform a Monte Carlo simulation study to evaluate the performance of the proposed tests. For comparison purposes, besides the proposed tests and the usual score and gradient tests, we also considered in the Monte Carlo experiment the improved score test (\citealt{n2010}), the modified score tests proposed by \citealt{k1996} and \citealt{cfc1998} and bootstrap-based tests. Our simulation results show that the improved gradient test proposed in this paper is an interesting alternative to the classic large-sample tests, delivering an accurate inference, mainly when dealing with small datasets. It is important to highlight that were not found any simulation study in the literature drawing a comparison between the performance of the proposed tests in the considered class of models, so this paper fills this gap.

The remainder of this paper is organized as follows. In Section $2$ we present the class of HSNLM, explaining inferential aspects. In Section $3$ we derive Bartlett and Bartlett-type corrections to improve the LR and gradient tests for testing varying dispersion in the model class of interest. We conduct a Monte Carlo study in order to evaluate and compare the performance of the proposed tests in Section $4$. An application to real data is presented in Section $5$. Some concluding remarks are given in Section $6$.

\section{Model specification}

Let $y$ be a random variable with symmetric distribution. Its density function is given by
\begin{equation}
\pi(y;\mu,\phi)=\frac{1}{\sqrt{\phi}}g(u), \ y, \ \mu \in \mathbb{R}, \ \phi>0,
\label{eq1}
\end{equation}
where $\mu$ is a location parameter, $\phi$ is a dispersion parameter, $u=(y-\mu)^2/\phi$, $g: \mathbb{R} \rightarrow [0,\infty)$ is the density generator (see, for example, \citealt{fkn1990}). We then denote $y \sim S(\mu, \phi,g).$ \citealt{CPG2005} presents the density generator function $g(\cdot)$ for some symmetric distributions.


Assume $y_1, \ldots, y_n$ being a random sample where each $y_\ell$ has a symmetric distribution (\ref{eq1}) with location parameter $\mu_\ell$ and dispersion parameter $\phi_\ell.$ Also, consider that the components of $\mat{\mu}=(\mu_1,\ldots, \mu_n)^\top$ and $\mat{\phi}=(\phi_1, \ldots, \phi_n)^\top$ vary across observations through nonlinear regression structures. The heteroscedastic symmetric nonlinear regression model $y_\ell \sim S(\mu_\ell, \phi_\ell, g), \ \ \ell=1,\ldots, n,$ proposed by \citealt{ccc2010} is defined by (\ref{eq1}) and by the systematic components for the mean vector response $\mat{\mu}$ and the dispersion parameter vector $\mat{\phi}$ described as follow $$\mu_{\ell}=f(\mat{x_\ell};\mat{\beta}) \ \ \mbox{and} \ \ \phi_\ell=h(\tau_\ell), $$ where $f(\cdot; \cdot)$ is a function possible nonlinear in the second argument which is continuous and differentiable in $\mat{\beta},$ where $\mat{\beta}=(\beta_1,\ldots, \beta_p)^\top$ ($p<n$ and $\mat{\beta} \in \mathbb{R}^p $) is a vector of unknown parameters to be estimated, and $\mat{x_\ell}=(x_{\ell 1},\ldots,x_{\ell m})^\top$ is a $m \times 1$ vector of known explanatory variables associated with the $\ell$th observation. Moreover, $h(\tau_\ell)$ is a known bijective continuously differentiable function of the scale linear predictor defined as $\tau_\ell=\mat{\omega_\ell}^\top \mat{\delta},$ where $\mat{\omega_\ell}=(1,\omega_{\ell 1},\ldots, \omega_{\ell k-1})^\top$ is a vector of explanatory variables which components are not necessarily different from $\mat{x_\ell},$ and $\mat{\delta}=(\delta_0,\ldots,\delta_{k-1})^\top$ $(\mat{\delta} \in \mathbb{R}^k)$ is a vector of unknown parameters. 

It is further assumed that if a value $\mat{\delta_0}$ of $\mat{\delta}$ exists, then $h(\mat{\omega_\ell}^\top\mat{\delta_0})=1$ for all $\ell$, therefore, $y_\ell's$ have constant dispersion if $\mat{\delta}=\mat{\delta_0}.$ The function $h(\cdot)$ should be a positive-value function and a possible choice is $h(\cdot) = \exp(\cdot),$ which is adopted in several papers (e.g., \citealt{cw1983}; \citealt{v1993}; \citealt{st1994}; \citealt{bc2005}). Furthermore, considering $h(\tau_\ell)=\exp(\tau_\ell)=\exp(\mat{\omega_\ell}^\top \mat{\delta}),$ it is not necessary impose any restriction on the components of $\mat{\omega}_\ell$ \cite[]{cw1983, lzx2009}. It is important to note that the meaning of heteroscedasticity we use in this work refers to varying dispersion, that is, when $\phi_1=\phi_2= \ldots=\phi_n$ we have a homoscedastic model; without this we have a heteroscedastic model.

Let $l(\mat{\theta})$ denote the total log-likelihood function for the parameter of vector $\mat{\theta}=(\mat{\beta}^\top, \mat{\delta}^\top)^\top$ given $y_1,\ldots, y_n.$ We have 
$l(\mat{\theta})=-\frac{1}{2}\sum_{\ell=1}^n\log(\phi_\ell)+ \sum_{\ell=1}^n t(z_\ell),$
with $t(z_\ell)=\log g(z_\ell^2)$ and $z_\ell=\sqrt{u_\ell}=\frac{(y_\ell - \mu_\ell)}{\sqrt{\phi_\ell}}.$ We assume that the function $l(\mat{\theta})$ is regular (\citealt[Chap 9]{ch1974}) with respect to all $\mat{\beta}$ and $\mat{\delta}$ derivatives up to fourth order. The total Fisher information matrix for $\mat{\theta}$ has a block diagonal structure, i.e., $\mat{K_\theta}=\mbox{diag}\{\mat{K_{\beta},\mat{K_{\delta}}} \},$ where $\mat{K_\beta}=-\alpha_{2,0}\mat{\tilde{X}}\mat{\Lambda}^{-1}\mat{\tilde{X}}$ and $\mat{K_\delta}=\mat{W^\top V W},$ with $\mat{\tilde{X}}=\partial \mat{\mu}/\partial \mat{\beta},$ $\mat{\Lambda}=\mbox{diag}\{1/\phi_1, \ldots, 1/\phi_n \},$ $\mat{W=\partial\mat{\tau}/\partial \mat{\delta}}$ and $\mat{V}=\mbox{diag}\{v_1,\ldots, v_n \},$ such that $v_\ell=((1-\alpha_{2,0}){h'}^2_\ell)/4\phi_\ell^2,$ where $h'=\partial \phi_\ell/\partial \tau_\ell$ and $\alpha_{r,s}=E\{t(z_\ell)^{(r)} z_\ell^s \}$ for $r,s \in \{1,2,3,4\}$ and $t(z_\ell)^{(k)}=\partial^k t(z_\ell)/\partial z_\ell^k,$ for $k=1,2,3,4$ and $\ell=1,\ldots,n.$ For some symmetric distributions, the quantities $\alpha_{r,s}$ are given in \citealt{ufc2008}. The parameters $\mat{\beta}$ and $\mat{\delta}$ are globally orthogonal, so their respective maximum likelihood estimators (MLEs), $\mat{\hat{\beta}}$ and $\mat{\hat{\delta}},$ are asymptotically independent. In order to obtain the MLEs $\mat{\hat{\beta}}$ and $\mat{\hat{\delta}}$ iteratively, the scoring method can be applied. This procedure is described in detail in \citealt{ccc2010}.

Our interest is to test heteroscedasticity in symmetric nonlinear regression models. The null and alternative hypothesis considered are, respectively, $H_0:\mat{\delta_1}=\mat{\delta_1^{(0)}}$ and $H_1:\mat{\delta_1}\neq\mat{\delta_1^{(0)}},$ where $\mat{\delta}$ is partitioned as $\mat{\delta}=(\delta_0, \mat{\delta_1}^\top)^\top,$ with $\delta_0$ a scalar and $\mat{\delta_1}=(\delta_1,\ldots,\delta_{k-1})^\top.$ Here, $\mat{\delta_1^{(0)}}$ is a fixed column vector of dimension $k-1$ such that $h(\mat{\omega_\ell^\top \delta_{1}^{(0)}})=1$ and $\delta_0$ and $\mat{\beta}$ are considered nuisance parameters. Actually, we are testing the dispersion parameters in HSNLM, considering that under the null hypothesis this model comes down to the symmetric nonlinear regression model. The partition previously considered for $\mat{\delta}$ induces the corresponding partitions: $\mat{W}=(\mat{W_0},\mat{W_1}),$ where $\mat{W_0}$ is an $n\times 1$ vector with all ones and $\mat{W_1}=\partial \mat{\tau}/\partial \delta_1,$
\[ \mat{K_\delta} = \left[ \begin{array}{cc}
  K_{\delta_0 \delta_0}   & \mat{K}_{\mat{\delta_0 \delta_1}} \vspace{0.5cm} \\
 \mat{K}_{\mat{\delta_1 \delta_0}} & \mat{K}_{\mat{\delta_1 \delta_1}}\end{array} \right ], \]
with $K_{\delta_0 \delta_0}=\mat{W_0^\top V W_0},$ $\mat{K}_{\mat{\delta_0  \delta_1}}^\top=\mat{K}_{\mat{\delta_1 \delta_0}}=\mat{W_1^\top V W_0}$ e $\mat{K}_{\mat{\delta_1 \delta_1}}=\mat{W_1^\top V W_1}.$ The likelihood ratio ($S_{LR}$), score ($S_r$) and gradient ($S_g$) statistics for testing $H_0$ can be expressed, respectively, as 
\begin{eqnarray*}
S_{LR} &=& 2\{l(\mat{\hat{\delta}_1}, \hat{\delta}_0, \mat{\hat{\beta}}) - l(\mat{\delta_1^{(0)}}, \tilde{\delta}_0, \mat{\tilde{\beta}}) \}, \\
S_r &=& \frac{1}{4}[\mat{W_1}\mat{\tilde{\Lambda}}(\mat{\tilde{S}}\mat{\tilde{F_1}}\mat{\tilde{u}}-\mat{\tilde{F_1}}\mat{\iota})]^\top(\mat{\tilde{R}}^\top\mat{\tilde{V}}\mat{\tilde{R}})^{-1}[\mat{W_1}\mat{\tilde{\Lambda}}(\mat{\tilde{S}}\mat{\tilde{F_1}}\mat{\tilde{u}}-\mat{\tilde{F_1}}\mat{\iota})] \ \mbox{e}\\
S_g &=& \frac{1}{2}[\mat{W_1}\mat{\tilde{\Lambda}}(\mat{\tilde{S}}\mat{\tilde{F_1}}\mat{\tilde{u}}-\mat{\tilde{F_1}}\mat{\iota})]^\top(\mat{\hat{\delta}_1}-\mat{\delta_1^{(0)}}),
\end{eqnarray*}
where $(\mat{\hat{\beta}},\hat{\delta_0},\mat{\hat{\delta_1}})$ and $(\mat{\tilde{\beta}},\tilde{\delta_0},\mat{\delta_1^{(0)}})$ are, respectively, the unrestricted and restricted (under $H_0$) MLEs of $(\mat{\beta},\delta_0,\mat{\delta_1}),$ $\mat{\iota}$ is an $n\times 1$ vector of ones and $\mat{R}=\mat{W_1}-\mat{W_0}\mat{C},$ with $\mat{C}=(\mat{W_0}^\top\mat{V}\mat{W_0})^{-1}(\mat{W_0}^{-1}\mat{V}\mat{W_1}).$ Under the null hypothesis, these statistics have an asymptotic $\chi^2_{k-1}$ distribution up to an error of order $n^{-1}.$

\section{Improved test inference}


In order to obtain a more accurate inference when dealing with small and moderate sized samples, some procedures based on second-order asymptotic theory have been developed in the literature. For the HSNLM, a Bartlett-type correction factor for the score statistic was derived by \citealt{n2010}. To provide another improved test statistics to test varying dispersion in the class of HSNLM, we will derive Bartlett and Bartlett-type correction factors for the LR and gradient statistics, respectively, considering the general procedures developed by \citealt{l1956} and \citealt{vfl2014}. The Bartlett and Bartlett-type correction factors are very general and need to be obtained for every model of interest, since they involve complex functions of the moments of log-likelihood derivatives up to fourth order. Details about the derivation of the Bartlett and Bartlett-type correction factors are given in Appendix A (Supplementary material). 

To test $H_0:\mat{\delta_1}=\mat{\delta_1^{(0)}}$ in HSNLM considering $h(\mat{\omega_l}^\top \mat{\delta})=\exp(\mat{\omega_l}^\top \mat{\delta}),$ i.e., the case of heteroscedasticity with multiplicative effects, the Bartlett-corrected $LR$ statistic is given by $$S_{{LR}^*}=\frac{S_{LR}}{1+ c/(k-1)}, $$ where $c= \epsilon(\mat{\delta})+\epsilon(\mat{\beta},\mat{\delta})-\epsilon(\delta_0) -\epsilon(\mat{\beta},\delta_0),$
{\small \begin{eqnarray*}
\epsilon(\mat{\delta})&=& N_1 tr \{\mat{Z_{\delta_{d}}}^{(2)}\} + N_2\mat{\iota}^\top \mat{Z_\delta}^{(3)} \mat{\iota} + N_3\mat{\iota}^\top \mat{\Lambda} \mat{Z_\delta}^{(3)}\mat{\iota} + N_4\mat{\iota}^\top \mat{\Lambda Z_\delta}^{(3)}\mat{\Lambda} \mat{\iota} \\
&+& N_5\mat{\iota}^\top \mat{Z_{\delta_{d}}}^{(2)}\mat{Z_\delta} \mat{\iota} + N_6 \mat{\iota}^\top  \mat{Z_{\delta_{d}}}^{(2)}\mat{Z_{\delta}} \mat{\Lambda} \mat{\iota} +(N_7+N_8)\mat{\iota}^\top \mat{\Lambda} \mat{Z_{\delta_{d}}}^{(2)}\mat{Z_{\delta}} \mat{\iota},\\
\\
\epsilon({\mat{\beta}})&=&- N_{15}tr\{ \mat{\Lambda}\mat{Z_{\beta_{d}}}\mat{Z_{\delta_{d}}} \} - (N_{10}+N_{12})\mat{\iota}^\top\mat{\Lambda} \mat{Z_{\beta_{d}}}\mat{Z_{\delta }}\mat{Z_{\delta_{d}}}\mat{\iota}\\
&+& N_{14}\mat{\iota}^\top\mat{\Lambda} \mat{Z_{\beta_{d}}}\mat{Z_{\delta}}\mat{Z_{\beta_{d}}}\mat{\Lambda} \mat{\iota} - (N_{11}+N_{13})\mat{\iota}^\top \mat{\Lambda} \mat{Z_{\beta_{d}}}\mat{Z_{\delta}}\mat{Z_{\delta_{d}}} \mat{\Lambda} \mat{\iota}\\
&+& N_9\mat{\iota}^\top\mat{\Lambda} \mat{Z_{\delta}}\mat{Z_{\beta}}^{(2)}\mat{\Lambda} \mat{\iota},\\
\\
\epsilon(\delta_{0})&=& N_1 tr \{\mat{Z_{\delta_{0d}}}^{(2)}\} + N_2\mat{\iota}^\top \mat{Z_{\delta_{0}}}^{(3)} \mat{\iota} + N_3\mat{\iota}^\top \mat{\Lambda} \mat{Z_{\delta_{0}}}^{(3)}\mat{\iota} + N_4\mat{\iota}^\top \mat{\Lambda Z_{\delta_{0}}}^{(3)}\mat{\Lambda} \mat{\iota} \\
&+& N_5\mat{\iota}^\top \mat{Z_{\delta_{0d}}}^{(2)}\mat{Z_{\delta_{0}}} \mat{\iota} + N_6 \mat{\iota}^\top  \mat{Z_{\delta_{0d}}}^{(2)}\mat{Z_{\delta_{0}}} \mat{\Lambda} \mat{\iota}\\ 
&+&(N_7+N_8)\mat{\iota}^\top \mat{\Lambda} \mat{Z_{\delta_{0d}}}^{(2)}\mat{Z_{\delta_0}} \mat{\iota} \ \ \ \ \ \mbox{and}\\
\\
\epsilon({\mat{\beta}})&=&- N_{15}tr\{ \mat{\Lambda}\mat{Z_{\beta_{d}}}\mat{Z_{\delta_{0d}}} \} - (N_{10}+N_{12})\mat{\iota}^\top\mat{\Lambda} \mat{Z_{\beta_{d}}}\mat{Z_{\delta_{0} }}\mat{Z_{\delta_{0d}}}\mat{\iota}\\
&+&N_{14}\mat{\iota}^\top\mat{\Lambda} \mat{Z_{\beta_{d}}}\mat{Z_{\delta_{0}}}\mat{Z_{\beta_{d}}}\mat{\Lambda} \mat{\iota} -(N_{11}+N_{13})\mat{\iota}^\top \mat{\Lambda} \mat{Z_{\beta_{d}}}\mat{Z_{\delta_{0}}}\mat{Z_{\delta_{0d}}} \mat{\Lambda} \mat{\iota}\\
&+& N_9\mat{\iota}^\top\mat{\Lambda} \mat{Z_{\delta}}\mat{Z_{\beta}}^{(2)}\mat{\Lambda} \mat{\iota},
\end{eqnarray*} }
\hspace{-.1cm}where $\mat{Z_\beta}=\mat{\tilde{X}}(\mat{\tilde{X}}^\top\mat{\Lambda}\mat{\tilde{X}})^{-1}\mat{\tilde{X}}^\top,$ $\mat{Z_\delta}=\mat{W}(\mat{W}^\top\mat{V}\mat{W})^{-1}\mat{W}^\top,$ $\mat{Z_{\delta_0}}=\mat{W_0}$ $(\mat{W_0}^\top \\ \mat{V}\mat{W_0})^{-1}\mat{W_0}^\top,$ $\mat{Z_\beta}^{(2)}=\mat{Z_\beta}\odot \mat{Z_\beta},$ $\mat{Z_\delta}^{(2)}=\mat{Z_\delta}\odot \mat{Z_\delta},$ $\mat{Z_{\delta_{0}}}^{(2)}=\mat{Z_{\delta_{0}}}\odot \mat{Z_{\delta_{0}}},$ $\mat{Z_\delta}^{(3)}=\mat{Z_\delta}^{(2)}\odot \mat{Z_\delta},$ $\mat{Z_{\delta_{0}}}^{(3)}=\mat{Z_{\delta_{0}}}^{(2)}\odot \mat{Z_{\delta_{0}}},$ $\odot$ denotes the Hadamard (elementwise) product of matrices, and $(\cdot)_d$ indicates that the off-diagonal elements of the matrix are set equal to zero. The elements $N_i, \ i=1,\ldots,15$ are scalars, given by
{\small \begin{eqnarray*}
N_1&=&\frac{1}{64} \left\{\alpha_{4,1}+6\alpha_{3,3}+17\alpha_{2,2}-1 \right\},\\
N_2&=&\frac{1}{64} \left\{ -245\alpha_{2,2}^2+ 496\alpha_{2,2}+3\alpha_{2,2}\alpha_{3,3} -3\alpha_{3,3}-251 \right\},\\
N_3&=&\frac{1}{64} \left\{-17\alpha_{2,2}^2 + 32\alpha_{2,2}-\alpha_{2,2}\alpha_{3,3} + \alpha_{3,3} -15 \right\}, \\
N_4&=&\frac{1}{384} \left\{\alpha_{3,3}^2-39\alpha_{2,2}^2+66\alpha_{2,2}-6\alpha_{2,2}\alpha_{3,3} + 10\alpha_{3,3}-23 \right\}, \\
N_5&=&\frac{1}{64} \left\{\alpha_{2,2}^2 -2\alpha_{2,2} + 1 \right\}, \ N_9=\frac{1}{8}\left\{\frac{\alpha_{3,1}^2}{\alpha_{2,0}^2}-4 \right\},\\
N_6&=&-\frac{5}{128}\left\{16\alpha_{2,2} -9\alpha_{2,2}^2 + \alpha_{3,3} - \alpha_{2,2}\alpha_{3,3}-7\right\},\\
N_7&=&\frac{1}{128}\left\{-43\alpha_{2,2}^2 + 80\alpha_{2,2} - \alpha_{2,2}\alpha_{3,3} + 3 \alpha_{3,3}-37 \right\},\\
N_8&=&\frac{1}{256}\left\{\alpha_{3,3}^2+\alpha_{2,2}^2+ 2\alpha_{2,2}+2\alpha_{3,3} + 2\alpha_{2,2}\alpha_{3,3}+1 \right\},\\
N_{10}&=&\frac{1}{16\alpha_{2,0}}\left\{(1-\alpha_{2,2})(24\alpha_{3,1}+3\alpha_{2,0}) \right\},\\
N_{11}&=&\frac{1}{64\alpha_{2,0}}\left\{(\alpha_{3,1}+2\alpha_{2,0})(9\alpha_{2,2}+\alpha_{3,3}-7) \right\}, \\
N_{12}&=&-\frac{5}{32\alpha_{2,0}}\left\{(\alpha_{3,1}+2\alpha_{2,0})(\alpha_{2,2}-1) \right\},\\
N_{13}&=&-\frac{1}{64\alpha_{2,0}}\left\{(\alpha_{3,1}+2\alpha_{2,0})(-7\alpha_{2,2}+\alpha_{3,3}+9) \right\},\\
N_{14}&=&\left\{ \frac{1}{16} \left[\frac{\alpha_{3,1}}{\alpha_{2,0}} \right]^2 + \frac{\alpha_{3,1}}{4\alpha_{2,0}} +\frac{1}{4} \right\}\ \textrm{and}\\ 
N_{15}&=&\frac{1}{8\alpha_{2,0}}\left\{\alpha_{4,2}+\alpha_{3,1}-4\alpha_{2,0} \right\}.
\end{eqnarray*} }

The improved gradient statistic is obtained by multiplying its original statistic by a polynomial in the original statistic itself. The corrected gradient statistic continues to have a chi-squared distribution under the null hypothesis but its asymptotic approximation error decreases from $n^{-1}$ to $n^{-2},$ providing a more accurate inference. To test $H_0:\mat{\delta_1}=\mat{\delta_1^{(0)}}$ in HSNLM when $h(\mat{\omega_\ell}^\top \mat{\delta})=\exp(\mat{\omega_\ell}^\top \mat{\delta}),$ the corrected gradient statistic is given by $$S_{g^*}=S_g \{1 - (c_g + b_gS_g + a_gS_g^2) \},$$ where $a_g=\frac{A_3^g}{12(k-1)((k-1)+2)((k-1)+4)}, \ b_g=\frac{A_2^g-2A_3^g}{12(k-1)((k-1)+2)},$$  $$\ c_g=\frac{A_1^g - A_2^g + A_3^g}{12(k-1)},$ with
{\small \begin{eqnarray*}
A_1^g &=& 12\alpha_{2,0}Q_2 \mat{\iota}^\top \mat{\Lambda}\mat{Z_{\beta}}^{(2)}\odot (\mat{Z_\delta}-\mat{Z_{\delta_0}})\mat{\Lambda}\mat{\iota} +3Q_2^2\mat{\iota}^\top\mat{\Lambda}\mat{Z_{\beta_d}}(\mat{Z_\delta}-\mat{Z_{\delta_0}})\mat{Z_{\beta_d}}\mat{\Lambda}\mat{\iota}\\
&+& 6Q_2^2\mat{\iota}^\top\mat{\Lambda}(\mat{Z_\delta} - \mat{Z_{\delta_0}})\odot \mat{Z_\beta}^{(2)}\mat{\Lambda}\mat{\iota} + 3Q_1Q_2\mat{\iota}^\top\mat{\Lambda}\mat{Z_{\beta_d}}(\mat{Z_\delta} - \mat{Z_{\delta_0}})\mat{Z_{\delta_{0d}}}\mat{\iota}\\
&+& 3Q_1Q_2\mat{\iota}^\top\mat{Z_{\delta_{0d}}}(\mat{Z_\delta}-\mat{Z_{\delta_0}})\mat{Z_{\beta_d}}\mat{\Lambda}\mat{\iota} + 3Q_1Q_2\mat{\iota}^\top(\mat{Z_\delta}-\mat{Z_{\delta_0}})_d(\mat{Z_\delta}-\mat{Z_{\delta_0}})\mat{Z_{\beta_d}}\mat{\Lambda}\mat{\iota}\\
&+& 6Q_1Q_2\mat{\iota}^\top(\mat{Z_\delta}-\mat{Z_{\delta_0}})_d\mat{Z_{\delta_0}}\mat{Z_{\beta_d}}\mat{\Lambda}\mat{\iota} + 3Q_1^2\mat{\iota}^\top(\mat{Z_\delta}-\mat{Z_{\delta_0}})_d(\mat{Z_\delta}-\mat{Z_{\delta_0}})\mat{Z_{\delta_{0d}}}\mat{\iota}\\
&+& 6Q_1^2\mat{\iota}^\top(\mat{Z_\delta}-\mat{Z_{\delta_0}})_d\mat{Z_{\delta_0}}\mat{Z_{\delta_{0d}}}\mat{\iota} + 3Q_1^2\mat{\iota}^\top\mat{Z_{\delta_{0d}}}(\mat{Z_\delta}-\mat{Z_{\delta_0}})\mat{Z_{\delta_{0d}}}\mat{\iota}\\
&+& 6Q_1^2\mat{\iota}^\top(\mat{Z_\delta}-\mat{Z_{\delta_0}})\odot\mat{Z_{\delta_0}}^{(2)}\mat{\iota} + 6Q_3tr\{ \mat{Z_{\delta_{0d}}}(\mat{Z_\delta}-\mat{Z_{\delta_0}})_d \}\\
&-&12Q_5tr\{\mat{\Lambda}(\mat{Z_\delta}-\mat{Z_{\delta_0}})_d \mat{Z_{\beta_d}} \} + 6Q_4 tr \{\mat{\Lambda}(\mat{Z_\delta}-\mat{Z_{\delta_0}})_d \mat{Z_{\beta_d}} \},
\end{eqnarray*}}
{\small \begin{eqnarray*}
A_2^g &=& -3Q_1Q_3\mat{\iota}^\top (\mat{Z_\delta}-\mat{Z_{\delta_0}})_d(\mat{Z_\delta}-\mat{Z_{\delta_0}})\mat{Z_{\beta_d}}\mat{\Lambda}\mat{\iota}\\
 &-&3Q_1^2\mat{\iota}^\top(\mat{Z_\delta}-\mat{Z_{\delta_0}})_d(\mat{Z_\delta}-\mat{Z_{\delta_0}}) \mat{Z_{\delta_{0d}}}\mat{\iota}\\
&-& 3Q^2_1\mat{\iota}^\top(\mat{Z_\delta}-\mat{Z_{\delta_0}})_d\mat{Z_{\delta_0}}(\mat{Z_\delta}-\mat{Z_{\delta_0}})_d\mat{\iota} - 6Q_1^2\mat{\iota}(\mat{Z_\delta}-\mat{Z_{\delta_0}})^{(2)}\odot\mat{Z_{\delta_0}}\mat{\iota}\\
&-& \frac{9}{4}Q_1^2\mat{\iota}^\top(\mat{Z_\delta}-\mat{Z_{\delta_0}})_d(\mat{Z_\delta}-\mat{Z_{\delta_0}})(\mat{Z_\delta}-\mat{Z_{\delta_0}})_d\mat{\iota} -\frac{3}{2}Q_1^2\mat{\iota}^\top(\mat{Z_\delta}-\mat{Z_{\delta_0}})^{(3)}\mat{\iota}\\
&-& 3Q_3tr\{(\mat{Z_\delta}-\mat{Z_{\delta_0}})_d^{(2)} \} \ \mbox{and}\\
\\
A_3^g &=& \frac{3}{4}Q_1^2\mat{\iota}^\top(\mat{Z_\delta}-\mat{Z_{\delta_0}})_d(\mat{Z_\delta}-\mat{Z_{\delta_0}})(\mat{Z_\delta}-\mat{Z_{\delta_0}})_d\mat{\iota}+\frac{1}{2}Q_1^2\mat{\iota}^\top(\mat{Z_\delta}-\mat{Z_{\delta_0}})^{(3)}\mat{\iota},
\end{eqnarray*}}
where $(\mat{Z_{\delta}}-\mat{Z_{\delta_0}})_d=\mat{Z_{\delta_{d}}}-\mat{Z_{\delta_{0d}}}$ and $Q_1, Q_2, Q_3, Q_4, Q_5$ are scalars given by:
{\footnotesize \begin{eqnarray*}
Q_1 &=& \frac{1}{8} \{1-3\alpha_{2,2} - \alpha_{3,3} \},\ \ \ \ Q_2 = -Q_5 = -\frac{1}{2} \{\alpha_{3,1} + 2\alpha_{2,0} \},\\
Q_3 &=& \frac{1}{16} \left \{ 7\alpha_{2,2} - 1 + 6\alpha_{3,3} + \alpha_{4,4} \right \},\ \ \mbox{e}\ \ Q_4 = \frac{1}{4} \left \{\alpha_{4,2} + 5\alpha_{3,1} + 4\alpha_{2,0}. \right \}
\end{eqnarray*}}

The correction factors which improve the LR and gradient statistics are not easy to interpret, although they involve only simple matrix operations and can be easily implemented in any programming environment which perform linear algebra operations, such as {\tt MAPLE}, {\tt Ox}, {\tt R}, etc. Also, they depend on the distribution in (\ref{eq1}) only through the $\alpha$'s and also depend on the number of nuisance parameters, the dimension of the hypothesis tested and the matrix $\mat{X}$ and $\mat{W}$ of covariates. Finally, all unknown parameters in the correction factors are replaced by their restricted MLEs.

\section{Numerical evidence}


The simulation experiments are based on the heteroscedastic symmetric nonlinear regression model
$$y_\ell=\beta_0+\exp\{ \beta_1 x_{\ell 1} \}+\sum_{s=2}^p \beta_s x_{s \ell} + \epsilon_\ell, \ \ell=1,\ldots,n,$$ 
where $\epsilon_\ell \sim S(0, \exp\{\mat{\omega_\ell}^\top \mat{\delta}\},g).$ The response variable was generated assuming that $\beta_0=\ldots=\beta_{p-1}=1,$ $\delta_0=0.1, \delta_2=0.3, \delta_3=0.5$ and $\delta_4=\delta_5=\delta_6=1$ and different values for $p$ and $k$ were considered. The covariates $x_1,\ldots, x_{p-1}$ and $\omega_1, \ldots, \omega_k$ were generated as random samples of the $U(0,1)$ distribution and were kept fixed throughout the simulations. The null hypothesis under test is $H_0:\delta_1=\ldots=\delta_{k-1}=0,$ i.e., $\exp\{\mat{\omega_\ell}^\top \mat{\delta}\}=\exp\{\delta_0\},$ that is, under $H_0$ we have constant dispersion. All results were obtained using 10,000 Monte Carlo replications. We also carried out an additional simulation study including bootstrap-based tests where we considered 500 bootstrap samples. The bootstrap sampling was performed parametrically under the null hypothesis. The simulation results are based on the Student-$t$ (with $\nu=5$) and power exponential (with $\kappa=0.3$) models. The following nominal levels and sample size were considered: $\alpha=1\%, 5\%$ and $10\%,$ and $n=20,\ 30,\ \mbox{and} \ 40,$ respectively. We shall report the null rejection rates of the tests based on the following statistics: the original likelihood ratio, score and gradient statistics ($S_{LR}, \ S_r, \ S_g$), their respective Bartlett and Bartlett-type corrected versions ($S_{{LR}^*},\ S_{{r}^*}, \ S_{g^*}$) and the monotonic versions of the corrected score statistic proposed by \citealt{k1996} and \citealt{cfc1998}  ($S_{{r_1}^*},\ S_{{r_2}^*}),$ respectively. The simulations were carried out using the {\tt Ox} matrix programming language (\citealt{d2006}). All entries are percentages. 




Tables \ref{tab1}-\ref{tab3} show results for different sample sizes while keeping fixed (varying) the number of nuisance (interest) parameters. The results clearly show that the LR test is notably liberal (i.e., it over-rejects the null hypotheses), especially when the number of interest parameters and nuisance parameters increase (the results varying the number of nuisance parameters are not shown to save space). It also can be noted that the gradient test behaves quite similar to the LR test, but is less size distorted, while the usual score  test performs much better than the other two uncorrected ones, although it is a bit liberal in a few cases. Considering $\alpha=1\%$ and $n=30$ for the Student-$t$ model (see Table \ref{tab1}), the null rejection rates for the $LR$ test are $3.2\% \ (k=3),$ $5.0\% \ (k=4)$ and $6.5\% \ (k=5),$ for the gradient test are $2.4\% \ (k=3),$ $3.7\% \ (k=4)$ and $5.3\% \ (k=5)$ and for the score test are $0.8\% \ (k=3),$ $1.0\% \ (k=4)$ and $1.0\% \ (k=5).$

The simulation results also showed that the corrected tests based on the $S_{{LR}^*},$ $S_{r^*}$ and $S_{g^*}$ statistics outperformed their uncorrected versions, independently of the sample size and the number of interest or nuisance parameters. Additionally, as shown in Tables \ref{tab1}-\ref{tab3}, the corrected versions of the LR and gradient tests are very sensitive to increasing the number of parameters in the model, whether they are interest or nuisance parameters. Otherwise, the corrected score test is not influenced by the increase in the number of parameters in the model and among the improved tests, the one based on the $S_{{r}^*}$ statistic presents the best performance, exhibiting null rejection rates very close to the nominal level in most cases. For example, considering the power exponential model (Table \ref{tab3}), if $k=3, n=20$ and $\alpha=10\%,$ the null rejection rate for the tests based on $S_{{LR}^*},$ $S_{{g}^*}$ and $S_{{r}^*}$ are, respectively, $16.7\%,$ $11.3\%$ and $10.3\%,$ while considering the same scenario with $k=4,$ the null rejection rates for the tests based on $S_{{LR}^*},$ $S_{{g}^*}$ and $S_{{r}^*}$ are, respectively, $21.0\%,$ $15.7\%$ and $9.6\%.$ Now considering the tests based on the monotonic versions of the corrected score statistics $S_{{r_1}^*}$ and $S_{{r_2}^*}$ proposed by Kakisawa (1996) and Cordeiro et al. (1998), the simulation results shows that the performance of the tests based on those statistics are very similar to the corrected score test, presenting the same null rejection rate in most cases. Finally, we can also observe that all corrected and uncorrected tests present null rejection rates very close to the corresponding nominal level as the sample size increases, as expected.

In order to evaluate the performance of the improved numerical tests, i.e., bootstrap-based tests, and compare it with the behavior of the uncorrected and analytical corrected tests, we developed a supplementary simulation study, presented in Table \ref{tab5}. The bootstrap versions of the of the LR, score and gradient tests, being $S_{LR}^{boot}$, $S_r^{boot}$ and $S_g^{boot}$ their respective test statistics, follow the steps described below. Considering the studied model under the null hypothesis, we generate $B$ bootstrap resamples ($y_1^*,\ldots,y_B^*$). In this step we replace the unknown parameter vector by its estimates obtained under the null hypothesis computed using the original sample $(y_1, \ldots, y_B).$ Then we calculate the statistic $S_i,$ $i=LR,s,g,$ for each pseudo sample $y_1^*,\ldots,y_B^*,$ denoting the resulting statistic by $S^{boot}_{i_b},$ $b=1,\ldots,B.$ It is worth noting that the resulting statistics $S_i^{boot}$  do not follows the $\chi^2$ distribution, and the tests based on these statistics perform as follows. We estimate the percentile $1-\alpha$ of $S^{boot}_{i_b}$ by $\hat{q}_{1-\alpha},$ such that $\#\{S^{boot}_{i_b} \leq \hat{q}_{1-\alpha}B \}/B=1-\alpha,$ where $\#$ denotes the set cardinality. One decide to reject the null hypothesis if $S_i> \hat{q}_{1-\alpha}.$ Another way is to state the decision rule based on the bootstrap $p$-value given by $p^*=\#\{ S^{boot}_{i_b} \geq S_i \}/B.$ As can be seen in Table \ref{tab5}, the bootstrap-based tests are less size distorted than the corresponding uncorrected tests. Also, for the LR and gradient tests, their bootstrap versions outperform the corrected ones. On the other hand, the bootstrap score test behaves, in general, similarly to the monotonic and non-monotonic corrected ones. Simulations considering different values of $n$ (not shown) exhibited a similar pattern. For example, considering $p=3,k=3,n=30$ and $\alpha=5\%$ (see Tables \ref{tab1}-\ref{tab3} for the non bootstrap-based tests), the null rejection rates are $11\% \ (S_{LR}),$ $6.4\% \ (S_{{LR}^*})$, $5.4\% \ (S_{LR}^{boot},S_g^{boot})$, $5.0\% \ (S_{r}),$ $4.8\% \ (S_{{r}^*}, S_{{r_1}^*}, S_{{r_2}^*}),$ $5.7\% \ (S_{r}^{boot}),$ $9.7\% \ (S_{g})$ and $5.4\% \ (S_{{g}^*})$ for the Student-$t$ model and $11.1\% \ (S_{LR}),$ $7.7\% \ (S_{{LR}^*}),$ $5.3\% \ (S_{LR}^{boot}),$ $5.2\% \ (S_r,S_{r}^{boot}, S_{g}^{boot})$ and $5.5\% \ (S_{{r}^*},S_{{r_1}^*},S_{{r_2}^*})$ for the power exponential model.

Completing our simulation study, we performed experiments to evaluate the power of the tests considering a grid of values for $\delta$. With the exception of the tests based on $S_{LR},$ $S_{{LR}^*},$ and $S_g$ which presented liberal behavior, all other tests studied in this paper were considered. The results are presented in Table \ref{tab7} and show that as $\delta$ increases the tests are more powerful, as expected. Also, the bootstrapped tests are less powerful than the others as the value considered for $\delta$ moves away from zero.

In summary, the simulation results presented in this section show that the $LR$ and gradient tests are considerably oversided (liberal) and the analytical Bartlett and Bartlett-type corrections for these tests are effective in reducing the size distortion. The score test is the best performing uncorrected test. Its (monotonic or not) corrected versions perform the same, being overall the best performing tests as along with all the bootstrapped tests.

\begin{table}[H]
\begin{center}
\caption{\small Null rejection rates $(\%)$ for $H_0: \delta_1=\ldots=\delta_k=0$ with $p=3;$ $t_5$ model.} 
\label{tab1}
\begin{tabular}{lclrrrrrrrrrrrr}
\hline
n&&$Stat$&&\multicolumn{3}{c}{$\alpha=10\%$}&&\multicolumn{3}{c}{$\alpha=5\%$}&&\multicolumn{3}{c}{$\alpha=1\%$}\\
&&&&\multicolumn{3}{c}{$k$}&&\multicolumn{3}{c}{$k$}&&\multicolumn{3}{c}{$k$}\\ \cline{5-7} \cline{9-11}\cline{13-15} \vspace{-0.3cm}\\
&&&&$3$&$4$&$5$&&$3$&$4$&$5$&&$3$&$4$&$5$\\ 
\hline
$20$&&\small{$S_{LR}$}&&$26.7$&$37.1$&$41.4$&&$17.7$&$26.6$&$30.2$&&$6.5$&$11.8$&$14.7$\\
&&\small{$S_{{LR}^*}$}&&$14.7$&$18.1$&$18.3$&&$8.2$&$10.7$&$10.7$&&$2.0$&$3.1$&$3.0$\\
&&\small{$S_r$}&&$11.0$&$11.3$&$12.2$&&$5.7$&$5.7$&$6.1$&&$0.9$&$1.1$&$1.0$\\
&&\small{$S_{r^*}$}&&$10.0$&$10.3$&$11.2$&&$5.5$&$5.2$&$5.8$&&$1.1$&$1.0$&$1.0$\\
&&\small{$S_{{r_1}^*}$}&&$10.1$&$10.3$&$11.2$&&$5.5$&$5.2$&$5.8$&&$1.1$&$1.0$&$1.0$\\
&&\small{$S_{{r_2}^*}$}&&$10.1$&$10.3$&$11.2$&&$5.5$&$5.2$&$5.8$&&$1.1$&$1.0$&$1.0$\\
&&\small{$S_g$}&&$24.1$&$33.1$&$36.5$&&$15.2$&$22.6$&$26.1$&&$5.3$&$9.2$&$11.4$\\
&&\small{$S_{g^*}$}&&$10.6$&$16.9$&$18.7$&&$5.6$&$10.2$&$11.5$&&$1.6$&$3.2$&$4.3$\vspace{0.2cm}\\

$30$&&\small{$S_{LR}$}&&$18.4$&$23.0$&$26.0$&&$11.0$&$14.1$&$17.3$&&$3.2$&$5.0$&$6.5$\\
&&\small{$S_{{LR}^*}$}&&$12.0$&$12.7$&$13.6$&&$6.4$&$6.7$&$7.8$&&$1.3$&$1.7$&$2.0$\\
&&\small{$S_r$}&&$10.4$&$10.2$&$11.0$&&$5.0$&$5.2$&$5.5$&&$0.8$&$1.0$&$1.0$\\
&&\small{$S_{r^*}$}&&$9.8$&$9.9$&$10.2$&&$4.8$&$5.1$&$5.1$&&$0.8$&$1.1$&$1.0$\\
&&\small{$S_{{r_1}^*}$}&&$9.8$&$9.9$&$10.2$&&$4.8$&$5.1$&$5.1$&&$0.8$&$1.1$&$1.0$\\
&&\small{$S_{{r_2}^*}$}&&$9.8$&$9.9$&$10.2$&&$4.8$&$5.1$&$5.1$&&$0.8$&$1.1$&$1.0$\\
&&\small{$S_g$}&&$17.1$&$20.1$&$24.0$&&$9.7$&$11.9$&$15.5$&&$2.4$&$3.7$&$5.3$\\
&&\small{$S_{g^*}$}&&$10.5$&$11.1$&$14.6$&&$5.4$&$6.0$&$8.3$&&$1.0$&$1.5$&$2.6$\vspace{0.2cm}\\

$40$&&\small{$S_{LR}$}&&$15.8$&$17.6$&$19.1$&&$9.9$&$10.7$&$11.1$&&$2.6$&$3.1$&$3.2$\\
&&\small{$S_{{LR}^*}$}&&$11.7$&$11.2$&$10.9$&&$6.1$&$6.0$&$5.2$&&$1.2$&$1.3$&$1.2$\\
&&\small{$S_r$}&&$10.8$&$10.5$&$10.6$&&$5.5$&$5.2$&$5.1$&&$0.9$&$0.9$&$0.9$\\
&&\small{$S_{r^*}$}&&$10.2$&$10.0$&$9.9$&&$5.3$&$5.1$&$4.8$&&$0.9$&$0.9$&$0.9$\\
&&\small{$S_{{r_1}^*}$}&&$10.2$&$10.1$&$10.0$&&$5.3$&$5.1$&$4.8$&&$0.9$&$0.9$&$0.9$\\
&&\small{$S_{{r_2}^*}$}&&$10.2$&$10.0$&$9.9$&&$5.3$&$5.1$&$4.8$&&$0.9$&$0.9$&$0.9$\\
&&\small{$S_g$}&&$15.2$&$16.6$&$17.8$&&$8.8$&$9.9$&$10.1$&&$2.1$&$2.5$&$2.5$\\
&&\small{$S_{g^*}$}&&$9.9$&$10.8$&$10.6$&&$4.9$&$5.5$&$5.4$&&$1.1$&$1.1$&$1.2$\\
\hline
\end{tabular} 
\end{center}
\end{table}

\begin{table}[H]
\begin{center}
\caption{\small Null rejection rates $(\%)$ for $H_0: \delta_1=\ldots=\delta_k=0$ with $p=3;$ power exponential $\kappa=0.3$ model.} 
\label{tab3}
\begin{tabular}{lclrrrrrrrrrrrr}
\hline
n&&$Stat$&&\multicolumn{3}{c}{$\alpha=10\%$}&&\multicolumn{3}{c}{$\alpha=5\%$}&&\multicolumn{3}{c}{$\alpha=1\%$}\\
&&&&\multicolumn{3}{c}{$k$}&&\multicolumn{3}{c}{$k$}&&\multicolumn{3}{c}{$k$}\\ \cline{5-7} \cline{9-11}\cline{13-15} \vspace{-0.3cm}\\
&&&&$3$&$4$&$5$&&$3$&$4$&$5$&&$3$&$4$&$5$\\
\hline
$20$&&\small{$S_{LR}$}&&$25.2$&$34.0$&$35.9$&&$16.3$&$23.7$&$25.3$&&$6.2$&$9.8$&$10.7$\\
&&\small{$S_{{LR}^*}$}&&$16.7$&$21.0$&$20.0$&&$9.6$&$12.5$&$12.0$&&$2.9$&$3.5$&$3.5$\\
&&\small{$S_r$}&&$9.7$&$9.9$&$10.3$&&$4.8$&$5.2$&$5.7$&&$0.9$&$1.7$&$1.6$\\
&&\small{$S_{r^*}$}&&$10.3$&$9.6$&$10.0$&&$5.4$&$3.6$&$4.5$&&$1.2$&$1.1$&$0.6$\\
&&\small{$S_{{r_1}^*}$}&&$10.3$&$10.1$&$10.1$&&$5.4$&$3.6$&$4.5$&&$1.2$&$1.1$&$0.6$\\
&&\small{$S_{{r_2}^*}$}&&$10.3$&$10.0$&$10.1$&&$5.4$&$4.3$&$4.7$&&$1.2$&$1.1$&$0.6$\\
&&\small{$S_g$}&&$23.0$&$30.1$&$31.5$&&$14.6$&$20.4$&$21.4$&&$5.1$&$7.6$&$8.7$\\
&&\small{$S_{g^*}$}&&$11.3$&$15.7$&$15.4$&&$6.2$&$8.5$&$9.0$&&$1.3$&$2.2$&$2.7$\vspace{0.2cm}\\

$30$&&\small{$S_{LR}$}&&$18.6$&$21.3$&$22.8$&&$11.1$&$13.1$&$14.2$&&$3.4$&$4.2$&$4.5$\\
&&\small{$S_{{LR}^*}$}&&$13.9$&$14.3$&$14.3$&&$7.7$&$8.0$&$7.5$&&$1.9$&$2.0$&$1.6$\\
&&\small{$S_r$}&&$9.9$&$10.3$&$10.1$&&$5.2$&$5.3$&$5.4$&&$1.4$&$1.4$&$1.3$\\
&&\small{$S_{r^*}$}&&$10.2$&$10.2$&$10.0$&&$5.5$&$4.5$&$4.8$&&$1.4$&$1.2$&$0.7$\\
&&\small{$S_{{r_1}^*}$}&&$10.2$&$10.3$&$10.0$&&$5.5$&$4.7$&$4.8$&&$1.4$&$1.2$&$0.7$\\
&&\small{$S_{{r_2}^*}$}&&$10.2$&$10.3$&$10.0$&&$5.5$&$4.7$&$4.8$&&$1.4$&$1.3$&$0.7$\\
&&\small{$S_g$}&&$18,0$&$20.2$&$21.5$&&$10.5$&$12.2$&$13.0$&&$3.2$&$3.5$&$3.7$\\
&&\small{$S_{g^*}$}&&$11.4$&$12.1$&$12.4$&&$6.0$&$6.4$&$6.2$&&$1.2$&$1.4$&$1.4$\vspace{0.2cm}\\

$40$&&\small{$S_{LR}$}&&$14.4$&$16.4$&$18.3$&&$8.2$&$9.3$&$10.9$&&$2.3$&$2.9$&$3.4$\\
&&\small{$S_{{LR}^*}$}&&$12.6$&$11.7$&$12.7$&&$6.6$&$6.4$&$7.3$&&$1.8$&$1.7$&$1.8$\\
&&\small{$S_r$}&&$10.4$&$9.3$&$10.1$&&$5.4$&$4.9$&$5.3$&&$1.2$&$1.1$&$1.4$\\
&&\small{$S_{r^*}$}&&$10.2$&$9.4$&$10.7$&&$4.8$&$4.8$&$5.7$&&$0.5$&$0.8$&$1.5$\\
&&\small{$S_{{r_1}^*}$}&&$10.3$&$9.4$&$10.7$&&$4.9$&$4.8$&$5.7$&&$0.8$&$0.9$&$1.5$\\
&&\small{$S_{{r_2}^*}$}&&$10.3$&$9.4$&$10.7$&&$4.9$&$4.8$&$5.7$&&$0.8$&$0.9$&$1.5$\\
&&\small{$S_g$}&&$13.9$&$15.5$&$17.3$&&$7.8$&$8.7$&$10.1$&&$2.1$&$2.5$&$2.8$\\
&&\small{$S_{g^*}$}&&$9.9$&$10.5$&$11.7$&&$5.0$&$5.6$&$6.4$&&$1.2$&$1.5$&$1.5$\\
\hline
\end{tabular} 
\end{center}
\end{table}

\begin{table}[H]
\begin{center}
\caption{\small Null rejection rates $(\%)$ for $H_0: \delta_1=\ldots=\delta_k=0$ with $p=3,$ $n=30;$ $t_5$ and power exponential $\kappa = 0.3$ models.} 
\label{tab5}
\begin{tabular}{lclrrrrrrrrrrrr}
\hline
Model&&$Stat$&&\multicolumn{3}{c}{$\alpha=10\%$}&&\multicolumn{3}{c}{$\alpha=5\%$}&&\multicolumn{3}{c}{$\alpha=1\%$}\\
&&&&\multicolumn{3}{c}{$k$}&&\multicolumn{3}{c}{$k$}&&\multicolumn{3}{c}{$k$}\\ \cline{5-7} \cline{9-11}\cline{13-15} \vspace{-0.3cm}\\
&&&&$3$&$4$&$5$&&$3$&$4$&$5$&&$3$&$4$&$5$\\ 
\hline
$t_5$ &&\small{$S_{LR}^{boot}$}&&$10.4$&$10.0$&$10.6$&&$5.4$&$4.8$&$5.5$&&$1.4$&$0.9$&$1.3$\\
&&\small{$S_r^{boot}$}&&$10.3$&$10.4$&$9.9$&&$5.7$&$5.3$&$5.0$&&$1.2$&$0.9$&$0.9$\\
&&\small{$S_g^{boot}$}&&$10.2$&$10.2$&$9.8$&&$5.4$&$5.0$&$5.2$&&$1.1$&$0.9$&$1.0$\vspace{0.2cm}\\

&&&&&&&&&&&&&\\
\small {Power exponential}&&&&&&&&&&&&&\\
&&\small{$S_{LR}^{boot}$}&&$10.2 $&$10.7$&$10.5 $&&$5.3 $&$5.4$&$5.5 $&&$1.3 $&$1.0$&$1.2 $\\
&&\small{$S_r^{boot}$}&&$10.3 $&$10.3$&$9.6 $&&$5.2 $&$5.1$&$5.0 $&&$0.9 $&$1.1$&$1.1 $\\
&&\small{$S_g^{boot}$}&&$10.3$&$10.3$&$9.7$&&$5.2$&$5.2$&$4.9$&&$0.9$&$1.0$&$0.9$\\
\hline
\end{tabular} 
\end{center}
\end{table}


\begin{table}[H]
\begin{center}
\caption{\small Non-null rejection rates $(\%)$ for $H_0: \delta_1=\ldots=\delta_3=\delta$ with $p=3,$ $n=30,$ $\alpha=10\%;$ $t_5$ and power exponential $\kappa=0.3$ models } 
\label{tab7}
\begin{tabular}{llcrrrrrrrr}
\hline
Model &Stat&&\multicolumn{8}{c}{$\delta$}\\ \cline{4-11}\vspace{-0.3cm}\\
&&&$0.5$&$1.0$&$1.5$&$2.0$&$2.5$&$3.0$&$3.5$&$4.0$ \\
\hline
$t_5$ &&$$&$$&$$&$$&$$&$$&$$&$$\\
&\small{$S_r$}&&$14.0$&$19.4$&$38.6$&$57.1$&$74.8$&$77.2$&$92.2$&$92.6$\\
&\small{$S_{r^*}$}&&$14.1$&$19.2$&$38.6$&$57.0$&$74.9$&$76.9$&$92.1$&$92.6$\\
&\small{$S_{{r_1}^*}$}&&$14.1$&$19.2$&$38.6$&$57.0$&$75.0$&$76.9$&$92.1$&$92.6$\\
&\small{$S_{{r_2}^*}$}&&$14.1$&$19.2$&$38.6$&$57.0$&$75.0$&$76.9$&$92.1$&$92.6$\\
&\small{$S_{g^*}$}&&$15.8$&$23.0$&$37.1$&$58.8$&$77.5$&$79.4$&$95.3$&$95.5$\\
&\small{$S_{LR}^{boot}$}&&$20.9$&$31.0$&$40.9$&$50.5$&$60.6$&$70.5$&$80.8$&$91.3$\\
&\small{$S_r^{boot}$}&&$20.3$&$30.3$&$41.0$&$50.9$&$60.8$&$70.7$&$80.7$&$90.6$\\
&\small{$S_g^{boot}$}&&$20.1$&$30.0$&$40.6$&$50.5$&$60.0$&$70.0$&$80.0$&$89.7$\vspace{0.2cm}\\

Power exponential&&$$&$$&$$&$$&$$&$$&$$&$$\\
&\small{$S_r$}&&$13.2$&$31.5$&$41.9$&$65.6$&$86.1$&$91.5$&$98.8$&$98.9$\\
&\small{$S_{r^*}$}&&$13.0$&$30.7$&$42.2$&$66.0$&$85.9$&$90.1$&$98.9$&$99.0$\\
&\small{$S_{{r_1}^*}$}&&$13.0$&$30.7$&$42.3$&$66.2$&$85.9$&$91.3$&$98.8$&$99.0$\\
&\small{$S_{{r_2}^*}$}&&$13.0$&$30.7$&$42.2$&$66.1$&$85.9$&$91.3$&$98.8$&$99.0$\\
&\small{$S_{g^*}$}&&$12.7$&$28.7$&$42.3$&$66.2$&$85.7$&$94.5$&$99.4$&$99.5$\\
&\small{$S_{LR}^{boot}$}&&$19.7$&$29.9$&$40.0$&$50.1$&$60.6$&$70.6$&$80.4$&$90.2$\\
&\small{$S_r^{boot}$}&&$20.1$&$30.5$&$40.9$&$51.1$&$61.3$&$71.6$&$81.6$&$91.5$\\
&\small{$S_g^{boot}$}&&$20.8$&$31.0$&$41.7$&$52.0$&$62.1$&$71.8$&$81.8$&$81.8$\\
\hline
\end{tabular} 
\end{center}
\end{table}

\section{Real data application}

In this section, we consider a dataset on weight of eye lenses of European rabbit in Australia (\textit{Oryctolagus Cuniculos}), $y,$ in $mg,$ and the age of the animal,  $x,$ in days, in a sample containing 71 observations. This dataset was analyzes by \citealt{w1998} (example 6.8) and \citealt{CPG2005} which showed some evidence of heteroscedasticity. The model considered in this article introduces a regression structure to model dispersion in the model proposed by \citealt{CPG2005}, being given by $$y_l=\exp \left( \beta_1 - \frac{\beta_2}{x_\ell + \beta_3} \right)e^{\epsilon_\ell}, $$ where $\epsilon_\ell \sim S(0, \exp\{\delta_1+\delta_2x_\ell\})$, $\ell=1,\ldots, 71.$ The main goal here is to test $H_0: \delta_2=0$ against $H_1: \delta_2 \neq 0.$ For this test, the observed values of the test statistics ($p-$values in parentheses) are: $S_{LR}=8.368 \ (0.004),$ $S_{LR^*}=8.348 \ (0.004),$ $S_r = 6.776 \ (0.009),$ $S_{r^*} = 6.678  \ (0.010),$ $S_{{r_1}^*}= 6.679 \ (0.010),$ $S_{{r_2}^*}=6.678 \ (0.010),$ $S_g=7.828 \ (0.005)$ and $S_{g^*}=7.430 \ (0.006).$ The $p-$value of the bootstrapped tests are: $S_{LR^*}^{boot}=0.005,$ $S_r^{boot}=0.012$ and $S_g^{boot}=0.008.$ Note that all tests that employ corrected and bootstrapped score statistics do not lead to rejection of the null hypothesis at the $1\%$ nominal level, while the tests that employ the other statistics lead to the opposite decision at the same nominal level. From our simulations, we concluded that the corrected tests outperform their uncorrected versions. Also, we noticed that the corrected score test, their monotonic and bootstrapped versions have the same behavior, presenting in most scenarios null rejection rates closer to the considered nominal level than the other tests, leading to a more reliable inference and being preferable.

\section{Concluding remarks}

In this paper we derive Bartlett and Bartlett-type corrections to improve hypothesis testing of the dispersion parameters for the class of HSNLM proposed by \citealt{ccc2010} and compare in simulation study the performance of the proposed tests with the score test, its Bartlett-type corrected version and the uncorrected LR and gradient tests. We also consider for the simulation study monotonic versions of the Bartlett-type corrected score test and bootstrapped tests.

The numerical evidence suggests that the usual LR and gradient tests have similar performance, being oversided, mainly if the sample size is small or even moderate. It is clear that the Bartlett and Bartlett-type corrections attenuate this tendency, but their effectiveness in correcting the size distortions of the tests are completely different. While the corrected LR test presents very distorted rejection rates, the corrected gradient test produces results comparable to those of the usual and (monotonic or not) Bartlett-type corrected score tests. Additionally, the corrected score test and the bootstrapped tests perform the best overall. An advantage of the analytically corrected tests in relation to the bootstrapped tests is that it does not demand much computational burden. Moreover, it is important to note that the corrected tests deliver more trustful inference than their uncorrected versions when dealing with small or even moderate sized sample. We hence recommend the use of the Bartlett-type corrected score and gradient or bootstrapped tests in applications.

\section*{Acknowledgements}
We are thankful for the financial support of CNPq and FACEPE. The research of Lourdes C. Montenegro was supported by CAPES-Brazil (Grant 6796/14-1).

\bibliographystyle{aabc}
\bibliography{mybibabv}

\begin{thebibliography}{34}
\providecommand{\natexlab}[1]{#1}

\bibitem[{Ara\'ujo et~al.(2020)Ara\'ujo,Cysneiros and Montenegro}]{acm2020}
\textsc{Ara\'ujo MC, Cysneiros AHMA and Montenegro LC}. 2020. Improved
  heteroskedasticity likelihood ratio tests in symmetric nonlinear regression
  models. Stat Pap 61: 167--188.

\bibitem[{Barroso and Cordeiro(2005)}]{bc2005}
\textsc{Barroso LP and Cordeiro GM}. 2005. Bartlett corretions in
  heteroskedastic $t$ regression models. Stat Probabil Lett 75: 86--96.

\bibitem[{Bartlett(1937)}]{b1937}
\textsc{Bartlett MS}. 1937. Properties of suficiency and statistical tests. P R
  Soc London 160: 268--282.

\bibitem[{Bayer and Cribari-Neto(2013)}]{bc2013}
\textsc{Bayer FM and Cribari-Neto F}. 2013. Bartlett corrections in beta
  regression models. J Stat Plan Infer 143: 531--547.

\bibitem[{Cook and Weisberg(1983)}]{cw1983}
\textsc{Cook D and Weisberg S}. 1983. Diagnostics for heteroscedasticity
  diagnostics in regression. Biometrika 70: 1--10.

\bibitem[{Cordeiro and Cribari-Neto(2014)}]{gc2014}
\textsc{Cordeiro GM and Cribari-Neto F}. 2014. An Introduction to Bartlett
  Correction and Bias Reduction. New York: Springer.

\bibitem[{Cordeiro and Ferrari(1991)}]{cf1991}
\textsc{Cordeiro GM and Ferrari SLP}. 1991. A modified score test statistic
  having chi-squared distribution to order $n^{-1}$. Biometrika 78: 573--582.

\bibitem[{Cordeiro et~al.(1998)Cordeiro,Ferrari and Cysneiros}]{cfc1998}
\textsc{Cordeiro GM, Ferrari SLP and Cysneiros AHMA}. 1998. A formula to
  improve score test statistics. J Stat Comput Sim 62: 123--136.

\bibitem[{Cox and Hinkley(1974)}]{ch1974}
\textsc{Cox DR and Hinkley DV}. 1974. Theoretical Statistics. London: Chapman
  and Hall.

\bibitem[{Cysneiros(2011)}]{n2010}
\textsc{Cysneiros AHMA}. 2011. Bartlett-type Correction in Heteroscedastic
  Symmetric Nonlinear Models. Valence: 26th International Workshop on
  Statistical Modelling.

\bibitem[{Cysneiros et~al.(2010)Cysneiros,Cordeiro and Cysneiros}]{ccc2010}
\textsc{Cysneiros FJA, Cordeiro GM and Cysneiros AHMA}. 2010. Corrected maximum
  likelihood estimators in heteroscedastic symmetric nonlinear models. J Stat
  Comput Sim 80: 451--461.

\bibitem[{Cysneiros et~al.(2005)Cysneiros,Paula and Galea}]{CPG2005}
\textsc{Cysneiros FJA, Paula GA and Galea M}. 2005. Modelos Simetricos
  Aplicados. Sao Paulo: ABE - XI Escola de Modelos de Regressao.

\bibitem[{Doornik(2006)}]{d2006}
\textsc{Doornik JA}. 2006. Ox: An Object-Oriented Matrix Programming Language 4
  ed.. London: Timberlake Consultants Ltd.

\bibitem[{Fang et~al.(1990)Fang,Kotz and Ng}]{fkn1990}
\textsc{Fang KT, Kotz S and Ng KW}. 1990. Symmetric Multivariate and Related
  Distributions. London: Chapman and Hall.

\bibitem[{Kakisawa(1996)}]{k1996}
\textsc{Kakisawa Y}. 1996. Higher order monotone bartlett-type adjustment for
  some multivariate test statistics. Biometrika 71: 233--244.

\bibitem[{Lange et~al.(1989)Lange,Little and Taylor}]{llt1989}
\textsc{Lange KL, Little RJA and Taylor JMG}. 1989. Robust statistical modeling
  using the $t$ distribution. J Am Stat Assoc 84: 881--896.

\bibitem[{Lawley(1956)}]{l1956}
\textsc{Lawley DN}. 1956. A general method for approximating to the
  distribution of the likelihood ratio criteria. Biometrika 71: 233--244.

\bibitem[{Lemonte(2011)}]{l2011}
\textsc{Lemonte AJ}. 2011. Local power of some asymptotic tests in expontial
  nonlinear regression models. J Stat Plan Infer 141: 1981--1989.

\bibitem[{Lemonte(2012)}]{l2012}
\textsc{Lemonte AJ}. 2012. Local power properties of some asymptotic tests in
  symmetric linear regression models. J Stat Plan Infer 142: 1178--1188.

\bibitem[{Lemonte(2013)}]{l2013}
\textsc{Lemonte AJ}. 2013. Nonnull asymptotic distributions of the lr, wald,
  score and gradient statistics in generalized linear models with dispersion
  covariates. Statistics 47: 1249--1265.

\bibitem[{Lemonte et~al.(2012)Lemonte,Cordeiro and Moreno}]{lcm2012}
\textsc{Lemonte AJ, Cordeiro GM and Moreno G}. 2012. Bartlett corretions in
  birnbaum-saunders nonlinear regression models. J Stat Comput Sim 82:
  927--935.

\bibitem[{Lemonte and Ferrari(2012a)}]{lf2012a}
\textsc{Lemonte AJ and Ferrari SLP}. 2012a. The local power of gradient test.
  Ann I Stat Math 64: 373--381.

\bibitem[{Lemonte and Ferrari(2012b)}]{lf2012b}
\textsc{Lemonte AJ and Ferrari SLP}. 2012b. Local power and size properties of
  the lr, wald, score and gradient tests in dispersion models. Stat Methodol 9:
  537--554.

\bibitem[{Lin et~al.(2009)Lin,Zhu and Xie}]{lzx2009}
\textsc{Lin JG, Zhu LX and Xie FG}. 2009. Heteroscedasticity diagnostics for
  $t$ linear regression models. Metrika 70: 59--77.

\bibitem[{Maior and Cysneiros(2018)}]{mc2018}
\textsc{Maior VQS and Cysneiros FJA}. 2018. Symarma: a new dynamic model for
  temporal data on conditional symmetric distribution. Stat Pap 59: 75--97.

\bibitem[{Medeiros and Ferrari(2017)}]{mf2017}
\textsc{Medeiros FMC and Ferrari SLP}. 2017. Small-sample testing inference in
  symmetric and log-symmetric linear regression models. Stat Neerl 71:
  200--224.

\bibitem[{Medeiros et~al.(2017)Medeiros,Ferrari and Lemonte}]{mfl2017}
\textsc{Medeiros FMC, Ferrari SLP and Lemonte AJ}. 2017. Improved inference in
  dispersion models. Appl Math Model 51: 317--328.

\bibitem[{Simonoff and Tsai(1994)}]{st1994}
\textsc{Simonoff JS and Tsai CH}. 1994. Use of modified profile likelihood for
  improved tests of constancy of variance in regression. Appl Stat-J Roy St C
  43: 357--370.

\bibitem[{Terrel(2002)}]{t2002}
\textsc{Terrel GR}. 2002. The gradient statistic. Comp Sci Stat 34: 206--215.

\bibitem[{Uribe-Opazo et~al.(2008)Uribe-Opazo,Ferrari and Cordeiro}]{ufc2008}
\textsc{Uribe-Opazo MA, Ferrari SLP and Cordeiro GM}. 2008. Improved score test
  in symmetric linear regression model. Commun Stat A-Theor 37: 261--276.

\bibitem[{Vargas et~al.(2013)Vargas,Ferrari and Lemonte}]{vfl2013}
\textsc{Vargas TM, Ferrari SLP and Lemonte AJ}. 2013. Gradient statistic:
  higher order asymptotics and bartlett-type corrections. Electron J Stat 7:
  43--61.

\bibitem[{Vargas et~al.(2014)Vargas,Ferrari and Lemonte}]{vfl2014}
\textsc{Vargas TM, Ferrari SLP and Lemonte AJ}. 2014. Improved likelihood
  inference in generalized linear models. Comput Stat Data An 74: 110--124.

\bibitem[{Verbyla(1993)}]{v1993}
\textsc{Verbyla AP}. 1993. Modelling variance heterogeneity: residual maximum
  likelihood and diagnostics. J Roy Stat Soc B Met 55: 509--521.

\bibitem[{Wei(1998)}]{w1998}
\textsc{Wei BL}. 1998. Exponential Family Nonlinear Models. Singapore:
  Springer-Verlag.

\end{thebibliography}

\end{document}